\documentclass{amsart}

\newtheorem{theorem}{Theorem}[section]
\newtheorem{lemma}[theorem]{Lemma}
\newtheorem{proposition}[theorem]{Proposition}
\newtheorem{corollary}[theorem]{Corollary}

\theoremstyle{definition}

\newtheorem{example}[theorem]{Example}

\theoremstyle{remark}
\newtheorem{remark}[theorem]{Remark}

\numberwithin{equation}{section}

\begin{document}

\title[Logarithmic Derivatives]{Logarithmic Derivatives of Solutions to
Linear Differential Equations}

\author{Christopher J. Hillar}
\address{Department of Mathematics, University of California, Berkeley, CA 94720.}
\email{chillar@math.berkeley.edu}
\thanks{This work is supported under a National
Science Foundation Graduate Research Fellowship.}



\begin{abstract}
Given an ordinary differential field $K$ of characteristic zero,
it is known that if $y$ and $1/y$ satisfy linear differential
equations with coefficients in $K$, then $y'/y$ is algebraic over
$K$.  We present a new short proof of this fact using Gr\"{o}bner
basis techniques and give a direct method for finding a polynomial
over $K$ that $y'/y$ satisfies. Moreover, we provide explicit
degree bounds and extend the result to fields with positive
characteristic. Finally, we give an application of our method to a
class of nonlinear differential equations.
\end{abstract}
\maketitle

\section{Introduction}
Using Gr\"{o}bner basis techniques, we provide new constructive
proofs of two theorems of Harris and Sibuya \cite{Har,Har2} (see
also, \cite{Sing,Spe} and \cite[Problem 6.60]{Stan}) that give
degree bounds and allow for several generalizations. Their results
were the following.

\begin{proposition}\rm  Let $N_{1},N_{2} > 1$ be positive integers
and let $K$ be a differential field of characteristic $0$.  Let
$F$ be a (differential) field extension of $K$ and suppose that
$L_{1}(Y)$ and $L_{2}(Y)$ are nonzero homogeneous linear
differential polynomials (of orders $N_{1}$ and $N_2$
respectively) with coefficients in $K$.  Further, suppose that one
of the following holds:
\begin{enumerate}
    \item $y \in F$ has $L_{1}(y) = L_{2}(1/y) = 0$.
    \item $N_{2} \leq q \in \mathbb Z_+$, and $y \in F$ has $L_{1}(y)
    = L_{2}(y^{q}) = 0$.
\end{enumerate}
Then, $Dy/y$ is algebraic over $K$.
\end{proposition}

In this article, we prove the following more refined result.
\begin{theorem}\label{mainthm}\rm  Let $N_{1},N_{2} > 1$ be positive integers
and let $K$ be a differential field of characteristic $p$.  Let
$F$ be a (differential) field extension of $K$.  Suppose that
$L_{1}(Y)$ and $L_{2}(Y)$ are nonzero homogeneous linear
differential polynomials (of orders $N_{1}$ and $N_2$
respectively) with coefficients in $K$.  Further, suppose that one
of the following holds:
\begin{enumerate}
    \item $p > N_{1}+N_{2}-2$ and $y \in F$ has $L_{1}(y) =
L_{2}(1/y) = 0$.
    \item $p > (q+1)(N_{1}-1)-1$, $N_{2} \leq q \in \mathbb Z_+$, and $y
\in F$ has $L_{1}(y) = L_{2}(y^{q}) = 0$.
\end{enumerate}
Then, $D^{j}y/y$ is algebraic over $K$ for all $j \geq 1$.
Moreover, the degree of the polynomial for $D^{j}y/y$
$(j=1,\ldots,N_1-1)$ in (1) is at most $\frac{{\left( {N_2  + N_1
- 2} \right)!}} {{\left( {N_2 - 1} \right)!}}$ and for (2), at
most $\frac{{\left( {q(N_{1}-1)+N_{1}-2} \right)!}} {{\left(
{q(N_{1}-1)- 1} \right)!}}$.  Additionally, if $p = 0$, then the
first condition in (1) and in (2) may be dropped.
\end{theorem}

This theorem is proved in section 2.  In section 3, we give
tighter degree bounds under the assumption that $K$ has
characteristic zero (or more generally, if $K$ is an infinite
perfect field), and finally in section 4, we describe how our
technique applies to certain nonlinear differential equations.

Recall that a polynomial $f \in K[x]$ is called \textit{separable}
if all of its roots are distinct, and a field $K$ is called
\textit{perfect} if every irreducible polynomial in $K[x]$ is
separable.  Examples of perfect fields include finite fields,
fields of characteristic zero, and, of course, algebraically
closed fields. It is interesting to note that there is a converse
to Theorem \ref{mainthm} for perfect fields $K$. Specifically, we
have the following result.

\begin{proposition}
Let $K$ be a perfect field.  If $y'/y$ is algebraic over $K$, then
both $y$ and $1/y$ satisfy linear differential equations over $K$.
\end{proposition}

\begin{proof}
Suppose that $K$ is perfect and $u = y'/y$ is algebraic over $K$.
Let $f(x) = x^m+a_{m-1}x^{m-1}+\cdots+a_0 \in K[x]$ be the monic,
irreducible polynomial for $u$ over $K$. Since $K$ is perfect, it
follows from basic field theory that gcd($f,\frac{{\partial
f}}{{\partial x}}) = 1$.  In particular, $\frac{{\partial
f}}{{\partial x}} \neq 0$.  Consider now,
\[0 = f(u)' = u'\left( {mu^{m - 1}  + \sum_{i = 1}^{m - 1}
{ia_i u^{i - 1} } } \right) + \sum_{i = 0}^{m - 1} {a_i 'u^i }.\]
Since $\frac{{\partial f}}{{\partial x}} = {mx^{m - 1} +
\sum\nolimits_{i = 1}^{m - 1} {ia_i x^{i - 1} } }$ is not the zero
polynomial, it follows that ${mu^{m - 1} + \sum\nolimits_{i =
1}^{m - 1} {ia_i u^{i - 1} } } \neq 0$ (by irreducibility of $f$).
Hence, $u' \in K(u)$ and the same holds for higher derivatives.

Next, notice that $(1/y)' = -y'/y^{2} = -u/y$ and an easy
induction gives us that $(1/y)^{(k)} = p_{k}(u,u',u'',\ldots)/y$,
in which $p_{k}$ is a polynomial (over $K$) in $u$ and its
derivatives (set $p_0 = 1$).  By above, the polynomials
$p_{k}(u,u',\ldots)$ lie in the field $K(u)$.  This implies that
they satisfy some (non-trivial) linear dependence relation,
\[ \sum\limits_{k = 0}^N {h_kp_k = 0}, \] in which $h_{k} \in
K$. Therefore,
\[ 0 = \sum\limits_{k = 0}^N {h_kp_k/y} = \sum\limits_{k =
0}^N {h_k (1/y)^{(k)} } \] as desired.  Performing a similar
examination on the derivatives of $y' = uy$ produces a linear
differential equation for $y$ over $K$, completing the proof.
\end{proof}

As an application of our main theorem, take $F$ to be the field of
complex meromorphic functions on $\mathbb C$ and $K = \mathbb Q$.
Then, the only $y$ such that both $y$ and $1/y$ satisfy linear
differential equations over $K$ are the functions, $y = ce^{ux},$
in which $u$ is an algebraic number of degree at most
min$\{N_1,N_2\}$ and $c \in \mathbb C \setminus \{0\}$.

This simple example shows that it is possible to produce a minimum
degree of min$\{N_1,N_2\}$ for $y'/y$; however, it is still an
open question of whether we can achieve a minimum degree close to
the bound given in Theorem \ref{mainthm} (or for characteristic
zero, Theorem \ref{charzerothm} below).

Theorem \ref{mainthm} can also be used to show that elements in a
differential field $F$ do not satisfy linear differential
equations over a subfield $K$, as the following example
demonstrates.

\begin{example}
\cite[Problem 6.59]{Stan} Let $K = \mathbb C(x)$ and $F = \mathbb
C((x))$. Then, sec$(x)$ does not satisfy a linear differential
equation over $K$. Indeed, cos$(x)'/$cos$(x) = -$tan$(x)$ is not
algebraic over $\mathbb C(x)$.
\end{example}

We begin with a lemma that is fundamental for the solution of this
problem, although its generality should be useful in many other
contexts. Let $K$ be a field and let $R = K[x_{1},\ldots,x_{n}]$
be a polynomial ring over $K$. Give $R$ a grading by assigning to
each $x_{i}$, a number $w(x_{i}) = w_{i} \in \mathbb N$, so that
\[ w\left( {\prod\limits_{i = 1}^n {x_i ^{v_i } } } \right) =
\sum\limits_{i = 1}^n {v_i w_i }.
\]

Then, we have the following more general version of a theorem of
Sperber \cite{Spe}.  Its proof can be found in \cite[Lemma
2.2.2]{Stu}.

\begin{lemma}\rm  Let $I$ be the ideal of $R$ generated by a
collection of polynomials and let $J$ be the (homogeneous) ideal
generated by the leading forms (with respect to the grading above)
of these polynomials. Then, dim$(J) \geq$ dim$(I)$.
\end{lemma}
From this we obtain the immediate corollary that we need.
\begin{corollary}\label{maincor}\rm  If $J$ is zero-dimensional, then so is $I$.
\end{corollary}

\section{The Main Theorems}

Before embarking on proofs of the theorems mentioned in the
introduction, we present a simple example to illustrate our
technique. Let $y_{1},y_{2},\ldots$ be variables. We will view
$y_{j} = D^{j}y/y$ as solutions to a system of polynomial
equations over $K[\{y_{j}\}_{j=1}^{\infty}]$. For example,
consider the system ($N_1 = 3, N_2 = 2$):
\begin{gather*}
y'''+a_2 y''+a_1 y' + a_0 y = 0,
\\
(1/y)''+b_1 (1/y)'+b_0 (1/y) = 0
\end{gather*}
in which $a_2,a_1,a_0,b_1,b_0 \in K$.  Dividing the first equation
by $y$ and expanding the second one gives us the more suggestive
equations:
\begin{gather*}
y_3 + a_2 y_2 + a_1 y_1 + a_0 = 0,
\\
(2y_{1}^{2}-y_{2})-b_1 y_{1}+b_0 = 0.
\end{gather*}
Also, differentiating the original equation for $1/y$ and
expanding, we have that
\begin{equation*}
-6y_1^{3}+6y_1y_2-y_3 + b_1(2y_{1}^{2}-y_{2})-y_1(b_{1}'+ b_0) +
b_0' = 0.
\end{equation*}
Thus, we may view $(y'/y,y''/y,y'''/y) = (y_{1},y_{2},y_3)$ as a
solution to a system of three polynomial equations in three
unknowns.

Let $w(y_i) = i$ define a grading of $K[y_1,y_2,y_3]$, and notice
that the system of leading forms, $\{y_3 = 0, 2y_{1}^{2}-y_{2} =
0, -6y_1^{3}+6y_1y_2 -y_3= 0\}$, has only the trivial solution
$(y_1,y_2,y_3) = (0,0,0)$.  In light of Corollary \ref{maincor},
it follows that the equations above define a zero-dimensional
variety, establishing Theorem \ref{mainthm} (1) for this example.

In general, we will construct a system of $N_1-1$ equations in
$N_1-1$ unknowns satisfied by the $y_i$.  These equations will
define a zero-dimensional variety, and thus, standard elimination
techniques (see, for instance, \cite{Cox}) give us a direct method
of computing, for each $i$, a nonzero polynomial (over $K$)
satisfied by $y_i$.

Let us first examine what happens when we compute $f_{n} =
D^{n}(1/y)$. Notice that
\begin{gather*}
f_{0} = 1/y
\\
f_{1} = -y^{-2}Dy = -y_{1}/y
\\
f_{2} = 2y^{-3}(Dy)^{2} -  y^{-2}D^{2}y = 2y_{1}^{2}/y - y_{2}/y
\\
f_{3} = -6y_{1}^{3}/y+6y_{1}y_{2}/y-y_{3}/y.
\end{gather*}
In general, these functions $f_{n}$ can be expressed in the form
$f_{n} = (1/y)p_{n}(y_{1},\ldots,y_{n})$ for polynomials $p_{n}
\in \mathbb Z[y_{1},\ldots,y_{n}]$.  More remarkably, with respect
to the grading $w(y_{i}) = i$, these $p_{n}$ are homogeneous of
degree $n$.  These facts are easily deduced from the following
lemma.

\begin{lemma}\label{recur1}\rm  Let $m \in \mathbb
Z_{+}$.  Then,
\[\frac{{p_m }} {{m!}} =  - \sum\limits_{j =
1}^{m - 1} {\frac{{p_{m - j} }} {{(m - j)!}}\frac{{y_j }} {{j!}}}-
\frac{{y_m }} {{m!}}.\]
\end{lemma}

\begin{proof}
Consider the following well-known identity (Leibniz' rule),
\[\sum\limits_{j = 0}^m {{m \choose j} \left( D^j h \right)\left(
D^{m-j}g \right)} = D^m(hg).\] Setting $h = y$ and $g = 1/y$, it
follows that \[ \sum\limits_{j = 0}^m {\frac{{D^j y}}
{{j!}}\frac{{D^{m - j} (1/y)}} {{(m - j)!}}}  = 0.\]

Multiplying the numerator and denominator by $y$ and rewriting
this expression gives us
\[ \frac{{p_m }} {{m!}} =  - \frac{{p_{m - 1} }} {{(m -
1)!}}\frac{{y_1 }} {{1!}} - \frac{{p_{m - 2} }} {{(m -
2)!}}\frac{{y_2 }} {{2!}} -  \cdots  - \frac{{p_1 }}
{{1!}}\frac{{y_{m - 1} }} {{(m - 1)!}} - \frac{{y_m }} {{m!}}.
\]
\end{proof}

We are now ready to prove Theorem \ref{mainthm} (1).

\begin{proof}[Proof of Theorem \ref{mainthm} (1)]
With $N_{1},N_{2}$ as in Theorem \ref{mainthm}, we suppose $N_{1}
= n$, $N_{2} = m$.  Dividing through by $y$ in the first
differential equation for $y$ gives us
\begin{equation}\label{lineq1} y_n
= - a_{n - 1} y_{n - 1}  -  \cdots  - a_1 y_1  - a_0 \ \ \ \ a_{i}
\in K
\end{equation}
while multiplying the second one for $1/y$ by $y$ produces the
equation
\[ p_m  + b_{m - 1} p_{m - 1}  +
\cdots  + b_0 = 0 \ \ \ \ b_{i} \in K.
\]
Differentiating the original linear differential equation for $y$,
$k$ times, we will arrive at linear equations $y_{n+k} =
L_{k}(y_{1},\ldots,y_{n-1})$ in terms (over $K$) of
$y_{1},\ldots,y_{n-1}$ like (\ref{lineq1}) above (by repeated
substitution of the previous linear equations).  If we also
differentiate the equation for $1/y$, $k$ times, we will produce
another equation for the variables $y_{i}$.  More formally, we
have that
\[ D^{m+k}(1/y) + D^{k}(b_{m - 1}D^{m-1}(1/y)) + \cdots  + D^{k}(b_0/y) =
0\] produces the equation (by Leibniz' rule)
\[ D^{m+k}(1/y) +
\sum\limits_{i = 0}^{m - 1} {\sum\limits_{j = 0}^k { {k \choose j}
\left( {D^j b_i } \right)\left( {D^{k - j + i} \left( {1/y}
\right)} \right)} } = 0.\] So finally (after multiplying through
by $y$), it follows that
\begin{equation}\label{maineqs1}
P_{m+k} := p_{m+k} + \sum\limits_{i = 0}^{m - 1} {\sum\limits_{j =
0}^k { {k \choose j} \left( {D^j b_i } \right) p_{k - j + i} } }
 = 0.
\end{equation}

It is clear that the leading homogeneous forms of the $P_{m+k}$
(with respect to the grading above) are $p_{m+k}$. Consider now
the ring homomorphism $\phi: K[\{y_{i}\}_{i=1}^{\infty}]
\rightarrow K[y_{1},\ldots,y_{n-1}]$ defined by sending $y_{j}
\mapsto 0$ for $j \geq n$ and $y_{j} \mapsto y_{j}$ for $j < n$.
Let $\tilde P_{m+k}$ denote the polynomials produced by
substituting the linear forms $L_{i}$ for the variables $y_{n+i}$
($i = 0,1,\ldots$) into the polynomials, $P_{m+k}$.  The leading
homogeneous forms of the $\tilde P_{m+k}$ will just be $\tilde
p_{m+k} := \phi(p_{m+k})$ because we are substituting linear
polynomials with strictly smaller degree (corresponding to the
grading).  In light of Corollary \ref{maincor}, it suffices to
show that the $n-1$ equations (in the $n-1$ variables),
\begin{equation}\label{system1}
\tilde p_{m} = 0, \ \tilde p_{m+1} = 0, \ \ldots, \ \tilde
p_{m+n-2} = 0,
\end{equation}
are only satisfied by the point $(0,\ldots,0)$ to prove the claim.

As for the degree bounds in the theorem, simply notice that each
polynomial, $\tilde p_{m+k}$, has total degree $m+k$ (coming from
the $y_{1}^{m+k}$ term), so that Bezout's theorem gives us the
bound of $\frac{{\left( {m  + n  - 2} \right)!}} {{\left( {m  - 1}
\right)!}}$.

Suppose that $(y_{1},\ldots,y_{n-1}) \neq (0,\ldots,0)$ is a zero
of the system in (\ref{system1}); we will derive a contradiction.
Let $r \in \{1,\ldots,n-1\}$ be the largest integer such that $y_r
\neq 0$, and choose $t \in \{0,\ldots,m-1\}$ maximal such that
$\tilde p_{m-t} = 0$, $\tilde p_{m-t+1} = 0$, $\ldots$, $\tilde
p_{m} = 0$.  If $t = m-1$, then $\tilde p_{1} = -y_1 = 0$, and so
the recurrence in Lemma \ref{recur1} and (2.3) give us that $y_{i}
= 0$ for $i \in \{1,\ldots,n-1\}$, a contradiction. Thus, $t \leq
m-2$. Using Lemma 2.1 with $\phi$ (and the maximality of $r$),
examine the following identity:
\[
\frac{{\tilde p_{m-t + r - 1} }} {{(m-t + r - 1)!}} =  -
\frac{{\tilde p_{m-t + r - 2} }} {{(m-t + r - 2)!}}\frac{{y_1 }}
{{1!}} -  \cdots - \frac{{\tilde p_{m - t} }} {{(m -
t)!}}\frac{{y_{r-1} }} {{(r-1)!}} - \frac{{\tilde p_{m - t-1} }}
{{(m - t-1)!}}\frac{{y_{r} }} {{r!}}.\]

From (\ref{system1}) and the property of $t$ above, it follows
that $\frac{{\tilde p_{m - t-1} }} {{(m - t-1)!}}\frac{{y_{r} }}
{{r!}} = 0$.  Thus, $y_{r} = 0$ or $\tilde p_{m - (t+1)} = 0$; the
first possibility contradicts $y_r \neq 0$, while the second
contradicts maximality of $t$.

This proves that the equations (\ref{system1}) define a
zero-dimensional variety.  Finally, we remark that our methods
still hold if the characteristic of $K$ is larger than the
denominator $(m+n-2)$ above.  This completes the proof of the
theorem.
\end{proof}

The proof for Theorem \ref{mainthm} (2) is similar to the one
above, however, the recurrences as in Lemma \ref{recur1} are
somewhat more complicated. Let $n \in \mathbb N$, $q \in \mathbb
Z_{+}$ and examine $f_{n,q} = D^{n}(y^{q})$. It turns out that
$f_{n,q}$ = $y^{q}p_{n,q}(y_{1},\ldots,y_{n})$ in which $p_{n,q}
\in \mathbb Z[y_{1},\ldots,y_{n}]$ is homogeneous of degree $n$
(with respect to the grading $w(y_i) = i$).  This follows in a
similar manner as before from the following lemma.

\begin{lemma}\label{recur2}\rm  Let $p_{n,1} = y_{n}$ for $n \in \mathbb N$
($y_{0} = 1$).
Then, for all $m \in \mathbb N,q > 1$,
\[p_{m,q} = y_m  + \sum\limits_{j =
0}^{m - 1} {m \choose j} y_j p_{m - j,q-1}  .\]
\end{lemma}

\begin{proof}
Use Leibniz' rule as in Lemma \ref{recur1} with $h = y^{q-1}$ and
$g = y$.
\end{proof}

The next lemma will be used in the proof of Theorem \ref{mainthm}
(2), and it follows from a straightforward induction on $a$ (using
Lemma \ref{recur2}).

\begin{lemma}\label{philemma}\rm  Let $\phi$ be
as in the proof of Theorem \ref{mainthm} (1) and $n \geq 2$. Then,
for all $a \in \mathbb Z_+$ and $b \in \mathbb N$, we have $\phi
\left( p_{(a+1)(n-1)+b,a} \right) = 0$.
\end{lemma}

We now prove Theorem \ref{mainthm} (2).

\begin{proof}[Proof of Theorem \ref{mainthm} (2)]
With $N_{1},N_{2}$ as in Theorem \ref{mainthm}, we suppose $N_{1}
= n$, $N_{2} = m \leq q$. As before, the first differential
equation for $y$ gives us
\begin{equation}\label{lineq2} y_n
= - a_{n - 1} y_{n - 1}  -  \cdots  - a_1 y_1  - a_0 \ \ \ \ a_{i}
\in K
\end{equation}
while the second one for $y^{q}$ (after dividing through by
$y^{q}$) produces the equation
\[ p_{m,q}  + b_{m - 1} p_{m - 1,q}  +
\cdots  + b_0 = 0 \ \ \ \ b_{i} \in K.
\]
Differentiating the original linear differential equation for $y$,
$k$ times, produces linear equations $y_{n+k} =
L_{k}(y_{1},\ldots,y_{n-1})$ in terms (over $K$) of
$y_{1},\ldots,y_{n-1}$ like (\ref{lineq2}) above.  If we also
differentiate the equation for $y^{q}$, $k$ times, we will arrive
at another equation for the variables $y_{i}$:
\begin{equation*}
P_{m+k,q} := p_{m+k,q} + \sum\limits_{i = 0}^{m - 1}
{\sum\limits_{j = 0}^k { {k \choose j} \left( {D^j b_i } \right)
p_{k - j + i,q} } }
 = 0.
\end{equation*}

It is clear that the leading homogeneous forms of the $P_{m+k,q}$
(with respect to the grading above) are $p_{m+k,q}$.  Let $\phi$
be as in the proof of Theorem \ref{mainthm} (1), and let $\tilde
P_{m+k,q}$ denote the polynomials produced by substituting the
linear forms $L_{i}$ for the variables $y_{n+i}$ ($i =
0,1,\ldots$) into the polynomials, $P_{m+k,q}$.  If $\tilde
p_{m+k,q} := \phi(p_{m+k,q}) \neq 0$, then the leading homogeneous
form of $\tilde P_{m+k,q}$ is $\tilde p_{m+k,q}$ because we are
substituting linear polynomials with strictly smaller degree
(corresponding to the grading).

Consider the following system of equations (recall that $q \geq m$
and $n \geq 2$),
\begin{equation}\label{system2}
\tilde p_{m,q} = 0, \ \tilde p_{m+1,q} = 0, \ \ldots, \ \tilde
p_{(q+1)(n-1)-1,q} = 0.
\end{equation}
We claim that $(0,\ldots,0)$ is the only solution to
(\ref{system2}).  Suppose, on the contrary, that
$(y_{1},\ldots,y_{n-1}) \neq (0,\ldots,0)$ is a solution to
(\ref{system2}), and let $r \in \{1,\ldots,n-1\}$ be the largest
integer such that $y_r \neq 0$. Also, choose $t \in
\{1,\ldots,q\}$ minimial such that
\begin{equation}\label{system2witht}
\tilde p_{tr,t} = 0, \ \tilde p_{tr+1,t} = 0, \ \ldots, \ \tilde
p_{(t+1)r-1,t} = 0.
\end{equation}
Clearly $t \neq 1$, as then $\tilde p_{r,1} = y_r = 0$, a
contradiction.  Applying Lemma \ref{recur2} with $\phi$ (and
maximality of $r$), examine the equation,
\begin{equation}\label{finalthm2eqn}
\tilde p_{(t+1)r-1,t} = \tilde p_{(t+1)r-1,t-1} +  \cdots +
{(t+1)r-1 \choose r} y_{r}\tilde p_{tr-1,t-1}.
\end{equation}
Using Lemma \ref{philemma} (with $a = t-1$) and the maximality of
$r$, we have $\tilde p_{tr+b,t-1} = 0$ for all $b \in \mathbb N$.
Consequently, (\ref{finalthm2eqn}) and (\ref{system2witht}) imply
that $\tilde p_{tr-1,t-1} = 0$.  Repeating this examination with
$\tilde p_{(t+1)r-2,t}, \tilde p_{(t+1)r-3,t},\ldots, \tilde
p_{tr,t}$ (in that order) in place of $\tilde p_{(t+1)r-1,t}$ on
the left-hand side of (\ref{finalthm2eqn}), it follows that
$\tilde p_{tr-i,t-1} = 0$ for $i = 1,\ldots,r$.  This, of course,
contradicts the minimality of $t$ and proves the claim.

It now follows from Corollary \ref{maincor} that the variety
determined by \[\left\{\tilde P_{m,q} = 0,\ldots,\tilde
P_{(q+1)(n-1)-1,q} = 0\right\}\] is zero-dimensional.  Therefore,
by an elementary dimension argument, there exists a subset of at
most $n-1$ of these equations that generates a zero-dimensional
ideal. The degree bounds in the theorem now follow by noticing
that each polynomial, $\tilde P_{m+k,q}$, has total degree at most
$m+k$, so that Bezout's theorem gives us the bound of
$\frac{{\left( {q(n-1)+n-2} \right)!}} {{\left( {q(n-1)- 1}
\right)!}}$.

Finally, we remark that the methods above still hold if the
characteristic of $K$ is larger than $(q+1)(n-1)-1$.
\end{proof}

\section{Tighter Degree Bounds}

It is possible to refine the degree bounds in our main theorem if
we restrict our attention to fields of characteristic zero (see
Remark \ref{perfectremark} concerning the perfect field case).

\begin{theorem}\label{charzerothm}\rm  Assuming the hypothesis as in Theorem
\ref{mainthm} with
$K$ having characteristic zero, the degree of the polynomial for
$D^{j}y/y$ $(j = 1,\ldots,N_1-1)$ over $K$ in (1) is at most ${N_2
+ N_1 - 2 \choose N_1 - 1}$ and for (2), at most
${q(N_{1}-1)+N_{1}-2 \choose N_1 - 1}$.
\end{theorem}

\begin{proof}
We prove the result for (1) as the other case is the same. Let
$N_{1} = n$, $N_{2} = m$ and set $\tilde P_{m+k} \in
K[y_{1},\ldots,y_{n-1}]$ ($k = 0,\ldots,n-2$) to be the
polynomials in (\ref{maineqs1}) after substitution of the linear
forms, $y_{n+i} = L_{i}(y_{1},\ldots,y_{n-1})$. Corresponding to
the grading $w(y_j) = j$, the weight of each monomial in $\tilde
P_{m+k}$ is less than or equal to $m+k$.  Let $\tilde S$ be the
set of all (affine) solutions with coordinates in $\overline{K}$
to the system $\{\tilde P_{m+k} = 0 \}^{n-2}_{k=0}$.  Our first
goal is to bound the cardinality of $\tilde S$ by ${m+n  - 2
\choose n - 1}$.

Suppose that $\{y_{i,1},\ldots,y_{i,s}\}$ is the list of all $s$
distinct $i$-th coordinates of members of $\tilde S$. Since $K$ is
infinite, there exists $k_i \in K$ such that $y_{i,j} \neq k_i$
for $j = 1,\ldots,s$.  Now, let $x_1,\ldots,x_{n-1}$ be variables
and consider the new polynomials $F_{m+k} \in
K[x_{1},\ldots,x_{n-1}]$ produced by the substitution $y_i =
x_{i}^{i} + k_i$ in the $\tilde P_{m+k}$.  As the $n-1$ equations
$\tilde P_{m+k} = 0$ define a zero-dimensional variety, so do the
$n-1$ equations $F_{m+k} = 0$.

Let $S$ denote the set of all solutions with coordinates in
$\overline{K}$ to the system $\{F_{m+k} = 0 \}^{n-2}_{k=0}$. Since
the total degree of each $F_{m+k}$ is just $m+k$, we have by
Bezout's theorem, \[|S| \leq \frac{(m+n-2)!}{(m-1)!} =
(n-1)!{m+n-2 \choose n-1}.\]

Consider the (set-theoretic) map $\psi: S \rightarrow \tilde S$
given by $(x_1,\ldots,x_{n-1}) \mapsto (x_{1} +
k_{1},\ldots,x_{n-1}^{n-1} + k_{n-1})$.  It is easy to see that
\begin{equation}
\sum\limits_{s \in \tilde S} {\left| {\psi ^{ - 1} (s)} \right|} =
| {S} |.
\end{equation}
Let $(y_1,\ldots,y_{n-1}) \in \tilde S$.  By our choice of $k_i$,
the polynomial $h_i(x_i) = x_i^i + k_i - y_i$ has precisely $i$
distinct zeroes.  These $i$ roots are distinct since
characteristic zero implies that gcd$(h_i,h_i') = 1$. Hence,
$\left| {\psi ^{ - 1} (s)} \right| \geq (n-1)!$ for all $s \in
\tilde S$, and so from \[|\tilde S|(n-1)! \leq  |{S}| \leq {m+n-2
\choose n-1}(n-1)!,\] we arrive at the desired bound on $|\tilde
S|$.

Finally, suppose that $g_i(x) \in K[x]$ is the irreducible
polynomial for $y_i$ over $K$.  Since $K$ has characteristic zero,
this polynomial has distinct roots.  Thus, there are deg($g_i$)
distinct embeddings $\sigma: K(y_i) \to \overline{K}$ that are the
identity on $K$.  Moreover, each of these homomorphisms extends to
an embedding $\tilde{\sigma}: \overline{K} \to \overline{K}$
\cite[p. 233]{lang}. In particular, the deg($g_i$) points,
$(\tilde{\sigma} y_1,\ldots,\tilde{\sigma} y_{n-1})$, are all
distinct elements of $\tilde S$.  Thus, from our bound on $|\tilde
S|$, we must have
\[\text{deg}(g_i) \leq |\tilde S| \leq {m+n  - 2 \choose n - 1}.\]
This completes the proof of the theorem.
\end{proof}

\begin{remark}\label{perfectremark}
Examining the proof of Theorem \ref{charzerothm}, it should be
clear that the same result holds when $K$ is an infinite perfect
field and the hypotheses of Theorem \ref{mainthm} are satisfied.
\end{remark}

We should also note that the proof above generalizes to bound the
number of distinct solutions to certain systems of equations.
Specifically, we have the following interesting fact.

\begin{theorem}\label{gennumrootsthm}\rm  Let $w(y_j) = j$ be the grading as
above and let $K$ be a field of characteristic zero.  Let $m \in
\mathbb Z_+$ and suppose that $\{F_{m+k}(y_1,\ldots,y_{n-1}) =
0\}_{k = 0}^{n-2}$ is a zero-dimensional system of polynomial
equations over $K$ such that each monomial in $F_{m+k}$ has weight
less than or equal to $m+k$. Then, this system will have at most
${m+n - 2 \choose n - 1}$ distinct solutions with coordinates in
$\overline{K}$.
\end{theorem}

In principle, the number of solutions for a generic system with
conditions as in Theorem \ref{gennumrootsthm} can be found by a
mixed volume computation and Bernstein's Theorem (see \cite{Cox},
for instance).  This approach, however, seems difficult to
implement.

\section{Applications to Nonlinear Differential Equations}

In the proof of Theorem \ref{mainthm}, it is clear that the
important attributes of the recursions as in (\ref{lineq1}) are
that they reduce the degree and are polynomial in nature.  In
particular, it was not necessary that they were linear. For
example, the system,
\begin{gather*}
yy'''+a(y')^2+by^2 = 0,
\\
(1/y)''+c(1/y)'+d(1/y) = 0
\end{gather*}
gives us the recurrence $y_3 +ay_1^2+b = 0$, which has $y_3$
expressible as a polynomial in $y_1,y_2$ with strictly smaller
weight. Repeated differentiation of this equation, preserves this
property.  In general, let $h \in K[z_1,\ldots,z_n]$ be a
homogeneous polynomial (with respect to total degree) such that
each monomial $z^{\alpha} = z_1^{\alpha_1} \cdots z_n^{\alpha_n}$
has
\[ \sum\limits_{i = 1}^n {(i - 1)\alpha _i }  < n.\]
If the hypothesis of Theorem \ref{mainthm} are weakened to allow
$y$ to satisfy an equation of the form, $D^n y =
h(y,Dy,\ldots,D^{(n-1)} y)$, then the proof applies without
change. A generalization along these lines was also considered by
Sperber in \cite{Spe}, however, the techniques developed here give
us degree bounds just as in Theorem \ref{mainthm} and also allow
for a positive characteristic.

\section{Acknowledgement}

We would like to thank Michael Singer and Bernd Sturmfels for
several interesting and useful discussions about this problem.


\end{document}